\documentclass[12pt]{article}
\usepackage{amsmath,amssymb,amsthm}

\title{Rayleigh's Stretched String}
\author{Mark B. Villarino\\
Depto.\ de Matem\'atica, Universidad de Costa Rica,\\
2060 San Jos\'e, Costa Rica}

\date{January 2, 2008}

\newtheorem{thm}{Theorem}
\newtheorem{cor}{Corollary}

\numberwithin{equation}{section}

\makeatletter
\def\section{\@startsection{section}{1}{\z@}{-3.5ex plus -1ex minus
			  -.2ex}{2.3ex plus .2ex}{\large\bf}}
\def\subsection{\@startsection{subsection}{2}{\z@}{-3.25ex plus -1ex
			  minus -.2ex}{1.5ex plus .2ex}{\normalsize\bf}}
\renewcommand{\@dotsep}{200} 
\makeatother

\renewcommand{\geq}{\geqslant}  
\renewcommand{\leq}{\leqslant}  

\begin{document}

\maketitle

\begin{abstract}
We  obtain a rigorous \emph{a priori} upper and lower bounds to the exact period of the celebrated \textsc{Rayleigh} stretched-string differential equation.\end{abstract}



\section{Introduction}

In his monumental treatise \emph{The Theory of Sound} (see pages 45-46 \cite{LR}) Lord \textsc{Rayleigh} proposes the following problem:
 \begin{quote}
\emph{A horizontal elastic steel wire having length $2L_0$ and spring constant $\sigma$ is stretched to length $2L$.  A body, having a mass $m$ much greater than the mass of the wire, is tied to the center of the wire and is put in motion on a vertical line through the center of the wire.  The only forces on the body are produced by tension in the wire (there is no gravity, no damping).  It is required to study the motion of the body.}
\end{quote}(We have used \textsc{Agnew}'s formulation \cite{ag}, p. 17-18.)  Let $y(t)\equiv y$ be the vertical displacement of the mass $m$ in time $t$, where $y$ is positive, $0$, or negative according as the body is above, or at, or below the wire in equilibrium or neutral position.

By the law of \textsc{Hooke}, the magnitude of the elastic force, or the \emph{tension}, $T(y)$, exerted by\emph{ each half }of the stretched string is equal to:\begin{eqnarray*}
T(y)& = &\sigma\cdot \frac{\rm{stretch}}{\rm{original \ length}} \\
 & = &\sigma\cdot \frac{\sqrt{L^2+y^2}-L_0}{L_0}  
\end{eqnarray*} and therefore the \emph{vertical component} of the \emph{total }force, exerted by \emph{both halves}, which is the only component acting to move the mass $m$ is given by\begin{equation}
\label{ }
\rm{vertical \  component \  of \  T(y)}=-2\sigma\cdot \frac{\sqrt{L^2+y^2}-L_0}{L_0}\frac{y}{\sqrt{L^2+y^2}},
\end{equation}where we have employed our \emph{sign} conventions.

But, by \textsc{Newton}'s second law of motion, that same force, with \emph{opposite} sign, is given by $m\frac{d^2y}{dt^2}$.

Therefore, equating the two forces, and bringing the expression for the elastic force on the same side as the force given by \textsc{Newton}'s second law, we find that the algebraic sum of the forces is zero and dividing by $m$, we conclude that  \emph{the equation of motion} of the mass $m$ is given by:
\begin{equation}
\label{RS}
\boxed{\frac{d^2y}{dt^2}+\Bigl(\frac{2\sigma}{m}\frac{\sqrt{L^2+y^2}-L_0}{L_0\sqrt{L^2+y^2}}\Bigr)y=0}
\end{equation}
This differential equation \eqref{RS} is nonlinear and very complicated, and we have never seen an \emph{exact} treatment of it in the literature.

\textsc{Rayleigh}, himself, simplifies the equation by the following physical reasoning;

\begin{quotation}``The tension of the string in the position of equilibrium depends on the amount of the stretching to which it has been subjected.  In any other position the tension is greater; but we limit ourselves to the case of vibrations so small that the additional stretching is a negligible fraction of the whole.  On this condition, \emph{the tension may be treated as a constant.}"\end{quotation}

(Our italics.)We just saw that the magnitude of the tension, $T(y)$, (or the force) from each half of the string on $m$ is given by
\begin{equation}
\label{T}
T(y)=\sigma\frac{\sqrt{L^2+y^2}-L_0}{L_0}
\end{equation}provided the stretching is well within elastic limits.  So, \textsc{Rayleigh} \emph{effectively takes }$y=0$ in \eqref{T} so that the equation \eqref{RS} becomes
\begin{equation}
\label{LRS}
\boxed{\frac{d^2y}{dt^2}+\frac{2\sigma}{m}\frac{L-L_0}{L_0L}y=0,}
\end{equation} or, using the notation \eqref{T} (with $y=0$, i.e., $T:=T(0)=\sigma \frac{L-L_0}{L_0}$), \begin{equation}
\label{LRS1}
\ddot{y}+\frac{2T}{mL}y=0
\end{equation}where the ``dots" refer to time derivatives.  If we assume that the initial displacement at time $t=0$ is $y_0$, then  the solution to \eqref{LRS1} is 
\begin{equation}
\label{as}
\boxed{y(t)=y_0\cos\left(\sqrt{\frac{2T}{mL}}t\right)}
\end{equation}which is \emph{simple harmonic motion} with period $\overline{P}$ given by:
\begin{equation}
\label{aP}
\overline{P}=\frac{2\pi}{\sqrt{\frac{2T}{mL}}}
\end{equation}

\textsc{Rayleigh} does not discuss \emph{how close} his approximate period $\overline{P}$ \eqref{aP} is to the \emph{exact} period $P$.  Indeed, we have been unable to to find any error analysis of \textsc{Rayleigh}'s solution in the literature.  Yet, such a famous and classical differential equation merits an investigation into the accuracy of its approximate solution.

Therefore we offer the following theorem to fill this gap.

\begin{thm}\label{bounds}
Let $P$ be the true period of oscilation of the mass $m$.  Then \textsc{Rayleigh}'s approximate period $\overline{P}$ \textbf{\emph{overestimates}} the true period, $P$.  Indeed, the following inequalities are valid:
\begin{equation}
\label{bounds}
\boxed{\frac{\overline{P}}{\sqrt{1+\frac{\sigma m}{2TL_0}\cdot y_0^2}}
\leq P\leq \overline{P}}
\end{equation}
\end{thm}

  Our lower bound shows the dependence of the period on the \emph{reciprocal}, of the initial displacement, $y_0$, and tends to the upper bound as $y_0$ tends to $0$.


\section {The Formula for the Period}


The differential equation \eqref{RS} does not have the variable $t$ nor the first derivative,  $\dot{y}$ ,appearing explicity.  Therefore, we can transform it into a separable equation for the \emph{velocity}, \begin{equation}
\label{v}
v:=\dot{y}
\end{equation}since
\begin{equation}
\label{yddot}
\ddot{y}=\frac{dv}{dt}=\frac{dv}{dy}\frac{dy}{dt}=v\frac{dv}{dy}
\end{equation}Substituting the right-hand side of \eqref{yddot} into \eqref{RS}, transposing the expression for the tension to the right-hand side, separating the variables and forming the indefinite integral of both sides we obtain\begin{equation}
\label{ }
\frac{v^2}{2}=-\frac{\sigma}{2m}\left(\frac{y^2}{2L_0}-\sqrt{L^2+y^2}\right)+C
\end{equation}where $C$ is the constant of integration.  At time $t=0$, the conditions are $y(0)=y_0$ and $v(0)=0$.  Therefore, evaluating the constant $C$, recalling that $v=\frac{dy}{dt}$,  taking the square root of both sides, and using a \emph{negative} sign because the displacement, $y(t)$, is a \emph{decreasing} function during the first quarter oscilation, we obtain:\begin{equation}
\label{ }
\frac{dy}{dt}=-\sqrt{\frac{2\sigma}{m}(y_0^2-y^2)\left(\frac{1}{L_0}-\frac{2}{\sqrt{L^2+y^2}+\sqrt{L^2+y_0^2}}\right).}
\end{equation}Multiplying the quarter-period by $4$, we obtain the following formula for the full period, $P$, of oscilation of the mass:

\begin{thm}
The \textbf{true period}, $P$, of oscillation of the mass $m$ is given by:
\begin{equation}
\label{P}
\boxed{P=4\sqrt{\frac{m}{2\sigma}}\int_{0}^{y_0}\frac{1}{\sqrt{y_0^2-y^2}\sqrt{\frac{1}{L_0}-\frac{2}{\sqrt{L^2+y^2}+\sqrt{L^2+y_0^2}}}}~dy.}
\end{equation}
\qed
\end{thm}

This integral is very complicated for exact computation.

Indeed, if we make the change of variable$$z^2:=L^2+y^2, \ \  z_0^2:=L^2+y_0^2$$ we obtain the formula\begin{equation}
\label{ }
P=4\sqrt{\frac{m}{2\sigma}}\int_{L}^{z_0}\frac{z}{\sqrt{\frac{-1}{2L_0}z^4+z^3+\left(\frac{l^2}{2L_0}-z_0\right)z^2-L^2z+L^2z_0+\frac{L^2}{2L_0}}}~dz.
\end{equation}   If we write\begin{equation}\frac{}{}
\label{ }
\sqrt{\frac{-1}{2L_0}z^4+z^3+\left(\frac{l^2}{2L_0}-z_0\right)z^2-L^2z+L^2z_0+\frac{L^2}{2L_0}}
\end{equation}as\begin{equation}
\label{ }
\sqrt{\frac{-1}{2L_0}}\sqrt{(z-a)(z-b)(z-c)(z-d)}
\end{equation}and \begin{equation}
\label{PP}
\mathbf{P}(z):=4\sqrt{\frac{m}{2\sigma}}\int\frac{z}{\sqrt{\frac{-1}{2L_0}}\sqrt{(z-a)(z-b)(z-c)(z-d)}}~dz
\end{equation}then\begin{equation}
\label{ }
\boxed{P=\mathbf{P}(z_0)-\mathbf{P}(L).}
\end{equation}The integral in \eqref{PP} is a very complicated algebraic combination of\emph{ an elliptic integral of the third kind and an elliptic integral of the first kind}, and we will not pursue this line of inquiry any further.

Instead, we turn our attention to the upper and lower bounds for the true period, $P$.


\section{Upper and Lower Bounds for the True Period}

\begin{proof}[Proof of the bounds on the true period]
We begin with the \emph{exact} formula \eqref{P} for the true period, $P$:\begin{equation}
\label{P}
\boxed{P=4\sqrt{\frac{m}{2\sigma}}\int_{0}^{y_0}\frac{1}{\sqrt{y_0^2-y^2}\sqrt{\frac{1}{L_0}-\frac{2}{\sqrt{L^2+y^2}+\sqrt{L^2+y_0^2}}}}~dy.}
\end{equation}

Taking $y=y_0=0$ in the sum ${\sqrt{L^2+y^2}+\sqrt{L^2+y_0^2}}$ only makes the integrand larger and we conclude that

\begin{align*}
\label{}
   P< &  4\sqrt{\frac{m}{2\sigma}}\int_{0}^{y_0}\frac{1}{\sqrt{y_0^2-y^2}\sqrt{\frac{1}{L_0}-\frac{1}{L}} }~dy\\
   =&  4\frac{\sqrt{\frac{m}{2\sigma}}}{\sqrt{\frac{T}{\sigma L}}}\frac{\pi}{2}\\
   =& \frac{2\pi}{\sqrt{\frac{2T}{mL}}}
\end{align*}and \emph{we have proved the validity of the \textbf{upper }bound}:\begin{equation}
\label{ub}
\boxed{P\leq \frac{2\pi}{\sqrt{\frac{2T}{mL}}}}.
\end{equation}

For the \emph{lower} bound, we proceed similarly, and we start by observing that if we take
 $y=y_0$ in the sum ${\sqrt{L^2+y^2}+\sqrt{L^2+y_0^2}}$ the integrand becomes smaller and we conclude that

\begin{align*}
\label{}
 P\geq &  4\sqrt{\frac{m}{2\sigma}}\int_{0}^{y_0}\frac{1}{\sqrt{y_0^2-y^2}\sqrt{\frac{1}{L_0}-\frac{1}{\sqrt{L^2+y_0^2}}} }~dy\\
   =&  4\sqrt{\frac{m}{2\sigma}}\frac{\pi}{2}\frac{1}{\sqrt{\frac{1}{L_0}-\frac{1}{\sqrt{L^2+y_0^2}}} }\\
   \geq&\frac{2\pi}{\sqrt{\frac{2\sigma}{m}}}\frac{1}{\sqrt{\frac{L + \frac{y_0^2}{2L}}{L_0\sqrt{L^2+y_0^2}}}}\\
   \geq&\frac{2\pi}{\sqrt{\frac{2T}{mL}+\frac{\sigma y_0^2}{LL_0}}}
\end{align*}where, twice, we have used the inequality$$\sqrt{L^2+y_0^2}\leq L+\frac{y_0^2}{2L}.$$

Therefore, \emph{we have proved the validity of the \textbf{lower} bound:}\begin{equation}
\label{lb}
\boxed{P\geq \frac{2\pi}{\sqrt{\frac{2T}{mL}+\frac{\sigma y_0^2}{LL_0}}}}.
\end{equation}This completes the proof of the theorem \ref{bounds}.
\end{proof}

\begin{cor}If $\overline{P}$ denotes \textsc{Rayleigh}'s approximate period \eqref{aP}, then, the \textbf{relative error} $$R:=\frac{P-\overline{P}}{P}$$ in the approximation $P\approx \overline{P}$, satisfies the inequality\begin{equation}
\label{ }
\boxed{-\frac{m}{4(L-L_0)}\cdot y_0^2\leq R\leq 0}
\end{equation}
\qed
\end{cor}

This shows that the lower bound to the relative error is proportional to the \emph{square} of the initial displacement and tends to zero with a ``quadratic" rate of convergence.

\subsubsection*{Acknowledgment}
 Support from the Vicerrector\'{\i}a de Investigaci\'on of the 
University of Costa Rica is acknowledged.


\end{document}